\documentclass[a4paper,12pt]{article}
\usepackage{amsmath}
\usepackage{amscd}
\usepackage{amsfonts}
\usepackage{amssymb}
\usepackage{amsthm}
\usepackage{graphicx}
\usepackage{color}
\usepackage[colorlinks=true,unicode]{hyperref}
\definecolor{darkblue}{rgb}{0.0,0.0,0.5}
\definecolor{darkred}{rgb}{0.5,0.0,0.0}
\hypersetup{
    colorlinks,
    citecolor=blue,
    filecolor=blue,
    linkcolor=red,
    urlcolor=blue
}

\newtheorem{Thm}{Theorem}[section]
\newtheorem{Lem}[Thm]{Lemma}
\newtheorem{Prop}[Thm]{Proposition}

\newtheorem{Con}[Thm]{Conjecture}

\theoremstyle{definition}
\newtheorem{Def}[Thm]{Definition}

\newtheorem{Rem}[Thm]{Remark}

\def\e {\varepsilon}

\def\r {\rho}

\def\a {\alpha}
\def\A {\mathcal{A}}

\def\b {\beta}
\def\B {\mathcal{B}}

\def\d {\delta}

\def\l {\lambda}

\def\o {\omega}
\def\om {\omega}
\def\s {\sigma}
\def\S {\Sigma}
\def\St {{\Sigma^2}}
\def\Stt {{(\Sigma^2)^2}}

\def\t {\tilde}
\def\ov {\overline}

\def\bbZ{\mathbb{Z}}

\def\bbR{\mathbb{R}}
\def\bbN{\mathbb{N}}

\def\mD{\mathfrak{D}}

\def\mlip{\mathop{\mathrm{Lip}}}
\def\f12{\frac{1}{2}}
\def\ds{dynamical system }
\def\dss{dynamical systems }

\def\st{such that }
\def\te{there exist }
\def\tes{there exists }

\def\diffeos {diffeomorphisms }

\def\homeo {homeomorphism }

\def\skpr {skew product }
\def\skprs {skew products }
\def\nbd {neighborhood }

\def\skpr {skew product }

\def\sa {statistical attractor }

\def\nbds {neighborhoods }

\def\Ast {A_{\textup{stat}}}

\newcommand{\K}{{\mathcal{K}}}

\def\ulim {\mathop{\ov\lim}}
\def\diam {\mathop{\mathrm{diam}}}

\title{Cascades of $\e$-invisibility}

\author{Yu.~Ilyashenko\thanks{Cornell University, US; Moscow State and Independent Universities, Steklov Math.
Institute, Moscow; \it{yulij@math.cornell.edu}} \and D.~Volk\thanks{Institute for Information Transmission Problems RAS, Independent University of Moscow, Laboratoire J.-V. Poncelet (UMI 2615); \it{dire.ulf@gmail.com}}}
\title{Cascades of $\e$-invisibility\thanks{Work partially supported by the grants NSF 0700973, RFBR 07-01-00017-a and carried out in Cornell University and Independent University of Moscow}}
\date{}

\begin{document}
\maketitle

\begin{center}
To Steve Smale who gifted so much new vision and inspiration to our community.
\end{center}

\begin{abstract}
We consider statistical attractors of locally typical dynamical systems and their ``$\e$-invisible'' subsets: parts of the attractors whose \nbds are visited by orbits with an average frequency of less than $\e \ll 1$. For extraordinarily small values of $\e$ (say, smaller than $2^{-10^6}$), an observer virtually never sees these parts when following a generic orbit.

A trivial reason for $\e$-invisibility in a generic \ds may be either a high Lipschitz constant ($\sim 1/\e$) of the mapping (i.e. it badly distorts the metric) or its close ($\sim \e$) proximity to the structurally unstable dynamical systems. However~\cite{bib:IN} provided a locally typical example of \dss with an $\e$-invisible set and a uniform moderate ($<100$) Lipschitz constant independent on $\e$. These \dss from~\cite{bib:IN} are also $|\log\e|^{-1}$-distant from structurally unstable dynamical systems (in the class $\mathcal S$ of skew
products). Recall that a property of \ds is locally typical if every close system possesses it as well. The invisibility property is thus $C^1$-robust.

We further develop the example of~\cite{bib:IN} to provide a better rate of invisibility while keeping the same radius of the ball in the space of skew products. Our construction is based on series of cascading dynamical systems. Each system incorporates the previous one and further boosts the invisibility rate. We give an explicit example of $C^1$-balls in the space of ``step'' skew products over the Bernoulli shift \st for each \ds from this ball a large portion of the \sa is invisible.
The systems have rate of invisibility $\varepsilon $ with
$\e = 2^{-n^k}$.

{\textbf keywords:} statistical attractor, structural stability, hyperbolicity, skew product, symbolic shift, invisibility

\end{abstract}

\section{Introduction}  \label{sec:intro}
An attractor of a dynamical system is a set of states to which the other states tend asymptotically. However, despite the simplicity of the idea, there are many
non-equivalent definitions of attractors. Formalizing the notion of attractor differently, one can obtain the
maximal and Milnor attractors~\cite{bib:Mil}, the non-wandering set and the Birkhoff center~\cite{bib:Katok}, as well as the statistical~\cite{bib:AAIS} attractor. Their definitions are not only formally different, but for certain (usually degenerate) \dss they describe different sets.

The notion of  the statistical attractor that is recalled in section~\ref{sec:1}, is one of the ways of describing what an observer will see if looking at a \ds for a long time. More precisely, this kind of attractor is the smallest closed set where orbits of generic points  concentrate in the sense of time averages: the proportion of time spent outside of any \nbd of the attractor tends to zero.

%

The paper is devoted to a new effect in the theory of dynamical
systems called \emph{invisibility of attractors.} The systems with
this property have large parts of attractors that can not be
observed in numerical experiments of any reasonable duration. On the
other hand, these systems have a moderate Lipschitz constant and
form a ball in the space of skew products of radius about~$\frac{C}{n^2}$. The parameter~$n$ characterizes the rate of invisibility~$\e$ that can be
%
%
made as small as $2^{-n}$.
Skew products from this ball are structurally stable.
We say, that an open set $R$ in the phase space is $\e
$-invisible and $\e $ is the rate of invisibility of the part of the
attractor that belongs to this set provided that there exists a set
of measure $\e^{\frac 1 2}$ such that any point outside this set
never visits $R$ under the $k$-th iterate of the map for $|\log \e |
< k < \e^{-\frac 1 2}$. In practice, take $n = 10^4, \e = 2^{-10^4}$, this
implies that an observer will never see an orbit that visits $R$
after $10^4$ of iterates. This effect was discovered in \cite{bib:IN}.

In the present paper, for any $n$ having the same meaning as above,
we construct an open set of skew products over the Bernoulli shift
that has a large part of attractor invisible with the rate of invisibility
$2^{-n^k}$ where $k$ is one third of the Hausdorff dimension of the
phase space. The natural parameter $\frac{1}{n^2}$ (up to a constant factor) is still
the radius of the ball in the space of skew products for which our
construction works.

When the results of~\cite{bib:IN} were presented to William Thurston, he
asked, whether it is possible to obtain the rate of invisibility as
a tower of exponents whose height grows with the dimension. Such a
rapid decay was not obtained, however, the double exponential decay
constructed above is a response to Thurston's challenge.

We construct our example as a sequence of dynamical systems of increasing dimension, the next
one is a skew product over the previous one. We refer to such strategy as ``cascading''.
While using this approach, one obtains the desired construction step by step, like ascending a staircase,
getting better rate of invisibility on the each step. 

\section{Main Theorem} \label{sec:1}

\begin{Def} \label{def:astat}
Let $(X, \mu)$ be a compact metric measure space and $F \colon X \to X$ be a homeomorphism. The \emph{\sa }of the \ds $(X, F)$ is the minimal closed set $\Ast \subset X$  such that for each open \nbd $U \supset \Ast$ almost every orbit spends almost all the time in $U$:
\begin{equation}    \label{eq:astat}
\lim_{n \to \infty} \frac1n \# \left\{ k \,|\, F^k(x) \in U, \, 0 \leq k < n \right\} = 1 \text{ \ for $\mu$-a.e. $x \in X$.}
\end{equation}
\end{Def}

\begin{Rem} \label{rem:astat}
Definition~\ref{def:astat} can be restated in the following equivalent way: the point $x \in X$ does \emph{not} belong to the \sa if and only if \tes an open \nbd $U \ni x$ \st almost every orbit visits $U$ with zero average frequency:
\begin{equation}    \label{eq:not_astat}
\lim_{n \to \infty} \frac1n \# \left\{ k \,|\, F^k(x) \in U, \, 0 \leq k < n \right\} = 0 \text{ \ for $\mu$-a.e. $x \in X$.}
\end{equation}
\end{Rem}

This Remark shows us why the \sa is always non-empty: $\Ast = X \setminus V$, where $V$ is the union of all the \nbds $U$ satisfying~\eqref{eq:not_astat}; and the compactness of $X$ implies that~\eqref{eq:not_astat} cannot hold for every open $U \subset X$.

\begin{Def} \label{def:invis}
An open set~$U$ is called \emph{$\e$-invisible} if almost every orbit visits $U$ with an average frequency $\e$ or less:
\begin{equation}    \label{eq:invis}
\ulim_{n \to \infty} \frac1n \# \left\{ k \,|\, F^k(x) \in U, \, 0 \leq k < n \right\} \leq \e \text{ \ for $\mu$-a.e. $x \in X$.}
\end{equation}
\end{Def}

\begin{Rem} \label{rem:invis}
Due to Remark~\ref{rem:astat}, each $U$ \st $U \cap \Ast = \emptyset$ is totally invisible ($\e = 0$).
\end{Rem}

Let~$I$ be the interval~$[-1, 2]$ and for any $k \ge 2$ consider a smooth embedding of the $k$-dimensional cube~$Q := I^k$ into a $k$-dimensional sphere $M := S^k$. Let $\mu$ be a smooth measure on $M$ \st $\mu|_Q$ is exactly the standard Lebesgue measure on a cube, $\mu(Q) = 3^k$. Let $\mD$ be the space of \diffeos $f \colon M \to M$ endowed with the $C^1$-metric. Denote by $\mD(L)$ the set of \diffeos $f\in \mD$ such that $\mlip f\le L, \mlip f^{-1} \le L$.

We denote by~$\Sigma^2$ the set of sequences of zeros and ones that are infinite both to the left and to the right:
$$
\Sigma^2 = \{ \om = (\ldots, \om_{-1}, \om_0, \om_1, \ldots)\,|\,\om_i\in\{0, 1\},\,i\in\bbZ\}.
$$
Let $\s \colon \S^2 \to \S^2$ be the Bernoulli shift:
$$
(\s \om)_i = \om_{i+1}.
$$
The set~$\S^2$ bears the standard metric
$$
d (\om^1, \om^2) = \begin{cases}
2^{-m},\text{ where } m = \min\{i \ge 0\,|\,\om^1_{-i} \neq \om^2_{-i} \text{ or } \om^1_i \neq \om^2_i \}, \text{ for } \om^1 \neq \om^2, \\
0,  \text{ for } \om^1 = \om^2,
\end{cases}
$$
and the standard Bernoulli measure~$\mu_\S$. Recall that~$\mu_\S$ is defined by its value on the cylinders
$$
\mu_\S \{\om\,|\, \om_{i_1} = j_1, \dots, \om_{i_m} = j_m \} = 2^{-m}.
$$
Note that this measure is invariant under $\s$. Now we consider the product of $k$ copies of Bernoulli shift:
$$
\ov \s \colon {(\Sigma^2)}^k \to {(\Sigma^2)}^k, \quad \ov \s  = \underbrace {\sigma \times \dots \times \sigma }_{k
\text{ times }}.
$$
We will never meet the elements of a single $\St$ later so we will use the letter~$\om$ for the elements of $(\S^2)^k$ and we will write just $\s$ instead $\ov \s$. For instance, $\om_0$ now is a vector of $k$ zeros and ones.

Now consider the metric measure space
$$
X := {(\S^2)}^k \times M,
$$
measure~$\mu_X$ on~$X$ being the Cartesian product of~$\mu_\S$ and $\mu$.
Note that the Hausdorff dimension of $X$ equals
$3k$; this justifies the description of $k$ in the abstract.
A step skew product~$F$ is defined as follows:
\begin{equation}       \label{eq:def_step_skew}
F \colon X \to X, \ (\om, x) \mapsto (\s \om, f_{\om_0}(x)), \ f_{\om_0}
\in \mD,
\end{equation}
Note that the fiber map~$f_{\o_0}$ depends only on the zero vector of the whole bi-infinite sequence $\o$. This dependence resembles step functions from which the term is borrowed. We denote by~$C^1_{k}$ the space of such step skew products equipped with the following metric:
\begin{equation}       \label{eq:def_metric}
d (F,G) = \max\limits_{\om_0} d_{C^1} (f_{\om_0}^{\pm 1}, g_{\om_0}^{\pm 1}).
\end{equation}
If each of $f_{\om_0}$ is within $\mD(L)$ for certain $L$ we will write $F \in C^1_{k} (L)$.

Also let $\pi , \pi_i$ be the projections
\begin{equation}   \label{eq:projk}
\pi \colon X \to M, \ (\om, x) \mapsto x, \ \text{ and } \
\pi_i \colon {(\S^2)}^k \times Q \to I, \ (\om, x_1, \dots , x_k) \mapsto x_i.
\end{equation}

\begin{Thm} \label{thm:1}
Consider any $n > 100$ and $k \ge 2$. Let $\nu = \frac1n$. There exists a ball $\mathcal B_p \subset C^1_{k} (L)$, $L < 2$, of radius $c \nu^2$ in sense of distance~\eqref{eq:def_metric}, the constant~$c$ independent on~$n$ and $k$. Each skew product $G \in \mathcal B_p$ has a statistical attractor whose large part belongs to an $\e$-invisible set $R$ for
$$
\e = 2^{-n^k},
$$
both the attractor and the invisibility with respect to~$\mu_X$.
In more detail, $\pi \Ast(G) \subset Q$, and
\begin{equation}    \label{eq:astat_detail}
Q^- := [5\nu, 1-5\nu]^k \subset \pi \Ast(G) \subset [-2 \nu, 1+2 \nu]^k =: Q^+;
\end{equation}
while the set
\begin{equation}    \label{eq:invis_detail}
R = \pi_k^{-1} \left(-2 \nu, \frac1{10}\right)
\end{equation}
is $\e$-invisible with the above $\e$.
\end{Thm}

\begin{Rem}
The crucial feature of this result is the independence of Lipschitz constant $L$ on $n$. It is easy to construct an example of $\e$-invisibility if we allow $L$ to depend on $n$. However, the \dss obtained this way will tend to degenerate systems as $n \to +\infty$. See~\cite{bib:IN} for more details on this subject.
\end{Rem}

\begin{figure}[hbt] \label{fig:thm1}
    \centering
    \includegraphics[width=300pt]{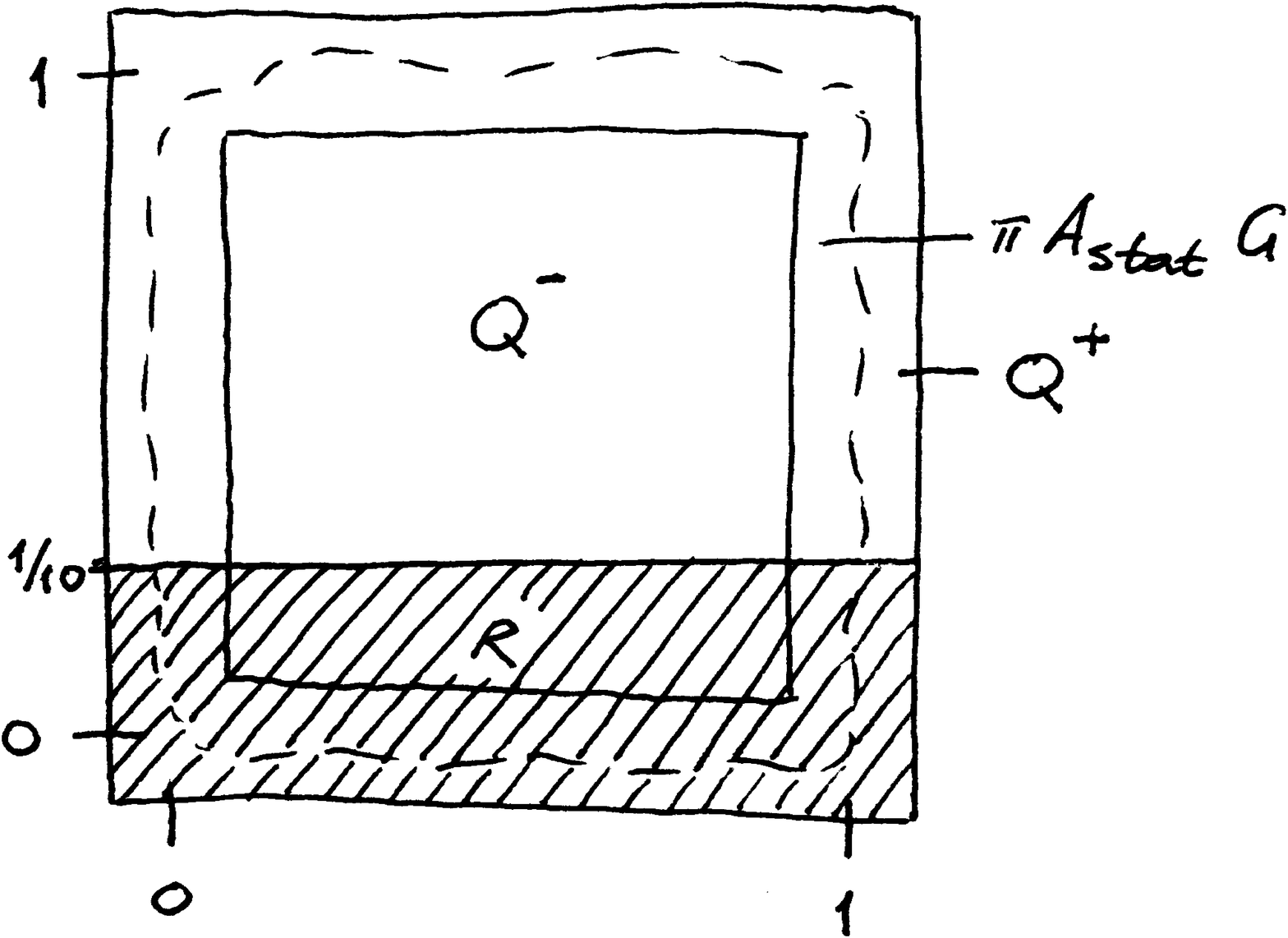}
    \caption{The \sa and the invisible set.}
\end{figure}

We also believe in the following smooth analogue of Theorem~\ref{thm:1}. Let $T \subset \bbR$ be a solid torus and $Y = T \times M$. Denote by $\mathcal{C}^1_k (L)$ the space of \diffeos~$\mathcal F$ of $Y$ such that $\mlip \mathcal F < L$, $\mlip \mathcal F^{-1} < L$, endowed with $C^1$-metric.

\begin{Con} \label{thm:smooth}
Consider any $n > 100$ and $k \ge 2$. Let $\nu = \frac1n$. There exists a ball $\mathcal B \subset \mathcal{C}^1_k (L)$, $L < 10$, $C \nu$-distant from structurally unstable diffeomorphisms. Each map $G \in \mathcal B$ has a large part of its \sa within $\e$-invisible set for
$$
\e = 2^{-n^k}
$$
in the same sense as in Theorem~\ref{thm:1}.
\end{Con}

First we give the detailed proof of Theorem~\ref{thm:1} for the case $k = 2$, as it is still simple enough and it contains all the techniques necessary for the general case. 
\section{Construction of the center of the ball} \label{sec:2}

In this section we explicitly construct the center~$F$ of the ball~$\mathcal B_p$ from Theorem~\ref{thm:1}. Recall that~$F$ is a step skew product. For $k = 2$ the fiber manifold~$M$ is a two-dimensional sphere $S^2$ and the base is $\Stt$. Hence we can define the skew product~$F$ of the form~\eqref{eq:def_step_skew} by fixing four \diffeos~$f_{00}, f_{01}, f_{10}, f_{11} \colon M \to M$.

\subsection{One-dimensional maps}

\begin{figure}[hbt]
    \centering
    \includegraphics[width=300pt]{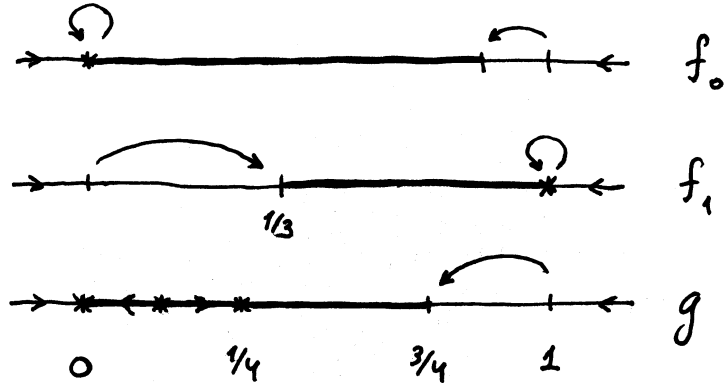}
    \caption{The maps~$f_0, f_1, g$. Bold segments are the images of $[0, 1]$.}
    \label{fig:factors}
\end{figure}

Recall that~$I = [-1, 2]$. Consider first one-dimensional orientation-preserving $C^1$-smooth maps $f_0, f_1, g \colon I \to I$, see~Figure~\ref{fig:factors}, with the following properties:

1) $f_0, f_1, g$ are \diffeos of $I$ onto its image that belongs to~$I$;

2) The map~$f_0$ has only one fixed point $x = 0$ and it is a weak attractor: the points of $I$ move towards $x = 0$ no more than by $\frac1{8n}$:
$$
\lim_{n \to \infty} f_0^n (x) = 0 \text{ and } \sup_{x \in I} | Id - f_0 | \le \frac1{8n};
$$

3) The map~$f_1$ has only one fixed point $x = 1$. We require this point to be a ``strong'' attractor with a multiplier independent of~$n$:
$$
\lim_{n \to \infty} f_1^n (x) = 1 \text{ and } f_1(I) \subset [\frac13, 1];
$$

4) The map~$g$ has $4n+1$ hyperbolic fixed points which are evenly spaced in the interval~$[0, \frac14]$ and the end points $x = 0$, $x = \frac14$ which are included into these $4n$ are attracting fixed points. We denote the distance between the adjacent fixed points of~$g$ by
\begin{equation}    \label{eq:def_h}
h := \frac{1}{16n};
\end{equation}

5) The Lipschitz constants of the \diffeos $f_0, f_1, g$ and their inverse maps~$f_0^{-1}, f_1^{-1}, g^{-1}$ are not greater than $\frac32$.

\begin{Prop}
The maps~$f_0, f_1, g$ with the properties~1~--~5 do exist.
\end{Prop}

\begin{proof}
We can define the maps~$f_0, f_1, g$ by the following formulas:
\begin{equation}    \label{eq:def_fg}
\begin{split}
& f_0(x) := \left(1 - \frac1{8n}\right) x, \\
& f_1(x) := 1 - \frac23 (1 - x), \\
& g(x) :=
\begin{cases}
\frac14 + \frac23 \left(x - \frac14\right), \ x \ge \frac14, \\
x - \frac{h}{3 \pi} \sin \frac{\pi x}{h}, \ 0 \le x \le \frac14, \\
\frac23 x, \ x \le 0.
\end{cases}
\end{split}
\end{equation}

It is easy to see that the properties~1~--~5 hold for these maps.
\end{proof}

\subsection{Two-dimensional maps}   \label{subsec:2-dim}

\begin{figure}[hbt]
    \centering
    \includegraphics[width=300pt]{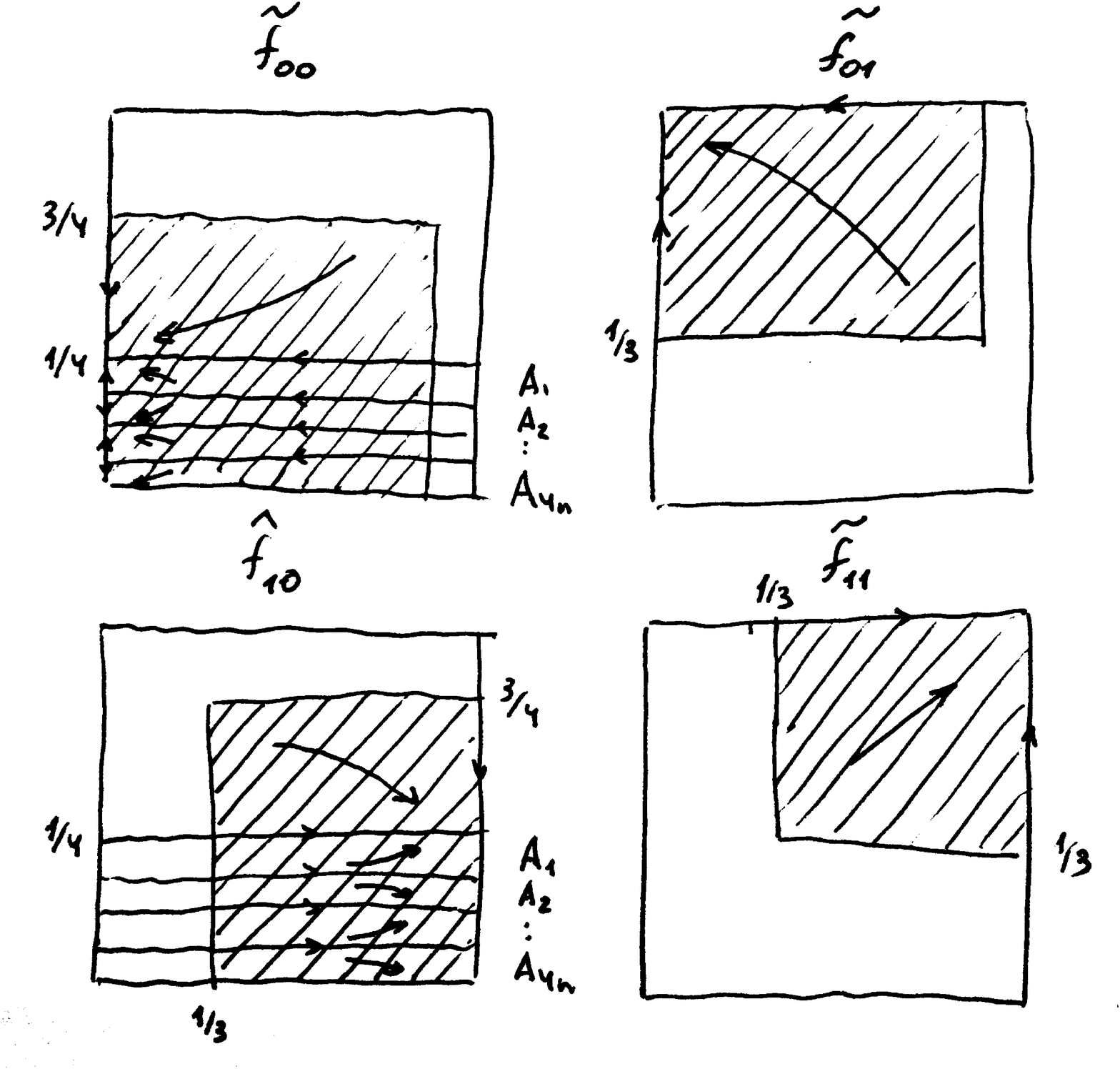}
    \caption{The maps~$\t f_{00}, \t f_{01}, \hat f_{10}, \t f_{11}$ on the square~$[0,1]^2$. The shaded regions are the images of $[0,1]^2$.}
    \label{fig:fij}
\end{figure}

Now we introduce the maps~$\t f_{00}, \t f_{01}, \t f_{10}, \t f_{11}$, which are \diffeos of the square~$Q^+$ onto its image. The final maps~$f_{ij}$ will be extensions of $\t f_{ij}$ onto the whole sphere~$M$:
$$
f_{ij}|_{Q^+} = \t f_{ij}.
$$
From now on, let~$f_0, f_1, g$ be the same as in~\eqref{eq:def_fg}. Also let $g_0(x) := g(x)$ and $g_1(x) := f_1(x)$. Then we define
\begin{equation}    \label{eq:def_fij}
\t f_{ij} := f_i \times g_i, \ ij \in \{0, 1\}^2, \ ij \neq 10, \ \text{ and } \hat f_{10} := f_1 \times g_0,
\end{equation}
see Figure~\ref{fig:fij}. The latter map is the prototype of the map~$\t f_{10}$.

\begin{figure}[hbt]
    \centering
    \includegraphics[width=300pt]{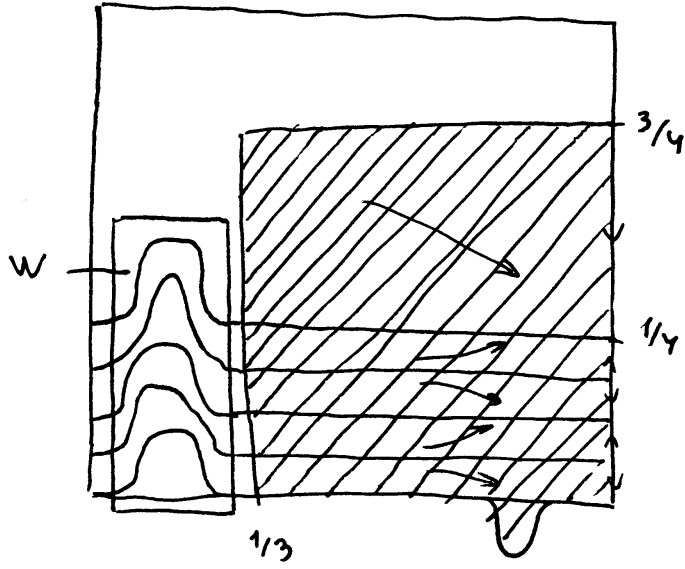}
    \caption{The map~$\t f_{10}$ on the square~$[0,1]^2$. The shaded region is the image of $[0,1]^2$.}
    \label{fig:f10}
\end{figure}

\begin{figure}[hbt]
    \centering
    \includegraphics[width=300pt]{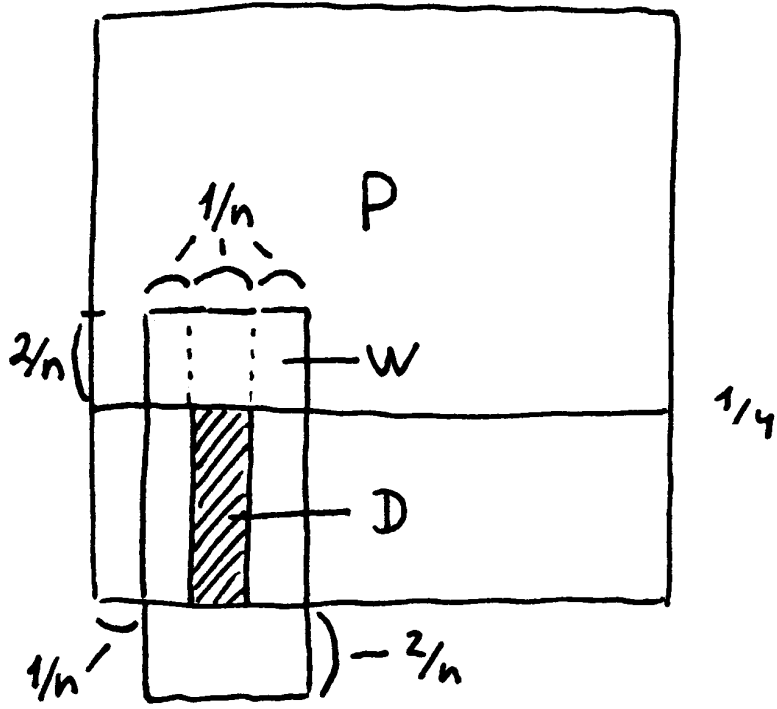}
    \caption{The fall-down region.}
    \label{fig:fall}
\end{figure}

Note that $\t f_{ij} (Q^+) \subset Q^+$, $ij \neq 10$. Points of~$Q^+$ uniformly
tend to some point of ${[0,1]}^2$ under the iterates of
any
of the maps~$\t f_{00}, \t f_{01}, \hat f_{10}, \t f_{11}$. Also note that the upper rectangle $P \subset [0,1]^2$,
\begin{equation}    \label{eq:def_p}
P := [0, 1] \times [\frac14, 1],
\end{equation}
is sent inside itself by every map~$\t f_{00}, \t f_{01}, \hat f_{10}, \t f_{11}$. Moreover, the maps~$\t f_{01}, \t f_{11}$ send the square~$[0,1]^2 \supset P$ inside~$P$. So the rectangle~$P$ in the fiber absorbs almost all the orbits from $\Stt \times Q^+$. The skew product~\eqref{eq:def_step_skew} requires a slight modification to destroy this property. We define
\begin{equation}    \label{eq:def_f10}
\t f_{10} (x, y) := (f_1 (x), g_0 (y) - \a (x) \b(y)),
\end{equation}
see Figure~\ref{fig:f10} and~\ref{fig:fall}, where $\a (x)$ and $\b (y)$ are defined as follows.

Let
\begin{equation}    \label{eq:def_dw}
D := \left[ \frac2n, \frac3n \right] \times \left[ 0, \frac14 \right], \ W := \left[ \frac1n, \frac4n \right] \times \left[ -\frac2n, \frac14 + \frac2n \right].
\end{equation}

Note that $D \subset W \subset Q^+$.

\begin{Prop}    \label{prop:ab_exist}
There exist functions~$\a \colon I \to I$, $\b \colon I \to [0,1]$ \st

1) $\a (x) \b (y) = \frac1{10n}$ for $(x, y) \in D$;

2) $\a (x) \b (y) = 0$ for $(x, y) \notin W$;

3) $\t f_{10}$ is a diffeomorphism, $\t f_{10} (Q^+) \subset Q^+$.

4) $\t f_{10}$ is Morse-Smale and $\frac1{10}$-distant from the structurally unstable diffeomorphisms.
\end{Prop}

\begin{proof}
Let $a < b < c < d$ and $b - a = d - c$. Define a $C^1$-smooth ``hat'' function~$\phi = \phi (\, \cdot \,; a, b, c, d) \colon I \to [0,1]$ \st
$$
\phi (x; a, b, c, d) := \begin{cases}
0, x \in I \setminus [a, d], \\
1, x \in [b, c], \\
\end{cases}
$$
and $|\phi'| \le \frac{1.01}{b - a}$. Take
$$
\a (x) := \phi \left(x; \frac1n, \frac2n, \frac3n, \frac4n\right),
$$
and
$$
\b(y) := \frac1{10n} \phi \left( y; -\frac2n, 0, \frac14, \frac14 + \frac2n\right).
$$
Properties~1 and~2 are obvious. Properties~3 and~4 will be checked in the proof of the subsequent Proposition. We only mention that
$$
|\a'\b| \le 1.01 \cdot n \frac1{10n} < \frac19, \quad |\a\b'| \le 1.01 \cdot \frac{n}{2} \frac1{10n} < \frac1{18},
$$
\begin{equation}    \label{eq:root}
\sqrt{|\a'\b|^2 + |\a\b'|^2} < \frac18.
\end{equation}
\end{proof}

The restrictions of the fiber maps~$f_{ij}$ to the square~$Q^+$ are now well defined.

\begin{Rem}
Now for each $ij \in \{0, 1\}^2$ we have $\t f_{ij} (Q^+) \subset Q^+$. Also note that each of the maps~$\t f_{ij}$ is a Morse-Smale diffeomorphism, which means they are structurally stable. We will use this feature in Sections~\ref{sec:6}-\ref{sec:9}.
\end{Rem}

\begin{Prop}    \label{prop:fij_exp_contr}
A) The maps~$\t f_{ij}$ are uniformly contracting on the rectangle~$P$, see~\eqref{eq:def_p};

B) The Lipschitz constants of the maps~$\t f_{ij}$ and of the inverse maps~$\t f_{ij}^{-1}$ are not greater than~$L = 1.85$.
\end{Prop}


The latter property is of little importance to this paper but it is essential for the proof of Theorem~\ref{thm:smooth} about the invisibility in \emph{smooth} case. A one-dimensional analogue of Theorem~\ref{thm:smooth} is proven in~\cite{bib:IN}; the proof involves an estimate similar to claim~B.

\begin{proof}
For $ij \neq 10$ the maps~$\t f_{ij}$ are Cartesian products, so the Proposition follows directly from the properties~1~--~5. In order to prove claim~A for $ij = 10$ we have to estimate the norm of $D \t f_{10}$, taken in the region~$P$. Denote $A := D \hat f_{10}$ and let
$$
B := \begin{pmatrix}
0 & 0 \\
\a'\b & \a\b' \\
\end{pmatrix}
$$
be the derivative of the map
$$
(x, y) \mapsto (0, \a(x)\b(y)).
$$
Note that~\eqref{eq:root} implies $\|B\| < \frac18$. The explicit calculation gives us within~$P$
$$
\|D \t f_{10}\| = \| A - B \| \le \|A\| + \|B\| < \frac23 + \frac18 < \frac56,
$$
which proves the claim~A about the contraction in~$P$.

The same argument gives us the following within~$Q^+$
\begin{equation}    \label{eq:185}
\|D \t f_{10}\| = \| A - B \| \le \|A\| + \|B\| \le \frac32 + \frac18 < 1.85,
\end{equation}
which provides us with the estimation of Lipschitz constant for the map~$\t f_{10}$. Let us also calculate the Lipschitz constant for the inverse map~$\t f_{10}^{-1}$:
$$
\|D \t f_{10}^{-1}\| = \| (A - B)^{-1} \| \le \|A^{-1}\| \cdot \| (Id - B A^{-1})^{-1} \| =
$$
$$
= \|A^{-1}\| \cdot \| Id + B A^{-1} + (B A^{-1})^2 + \dots \| \le \|A^{-1}\| \cdot \frac{1}{1 - \|B\| \cdot \|A^{-1}\|} \le \frac32 \cdot
\frac{1}{1 - \frac32 \cdot \frac18} < 1.85.
$$
Claim~B is verified too.

Now we prove claims~3 and~4 of Proposition~\ref{prop:ab_exist}. First, we show that the map~$\t f_{10}$ is a diffeomorphism, that is, globally invertible. It is enough to verify that $\| D \t f_{10} - Id \| < 1$:
$$
\| D \t f_{10} - Id \| = \| A - B - Id \| \le \| A - Id \| + \|B\| \le \frac12 + \frac18 < 1.
$$

Now we prove that $\t f_{10} (Q^+) \subset Q^+$. Due to~\eqref{eq:astat_detail} and~\eqref{eq:def_fg},
$$
\hat f_{10} (Q^+) \subset \left[ -\frac{4}{3n}; 1 + \frac{4}{3n}\right]^2.
$$
This implies that
$$
\t f_{10} (Q^+) \subset \left[ -\frac{4}{3n} - \| \a \b \|_{C^0}; 1 + \frac{4}{3n} + \| \a \b \|_{C^0}\right]^2.
$$
But $\| \a \b \|_{C^0} = \frac{1}{10n}$, so~$\t f_{10} (Q^+) \subset Q^+$.

Then we estimate the $C^1$-distance from $\t f_{10}^{\pm 1}$ to the structurally unstable (i.e. non-hyperbolic) \diffeos:
$$
\| D \t f_{10} - Id \| = \| A - B - Id \| \ge \| A - Id \| - \| B \| > \frac5{24} > \frac1{10};
$$
$$
\| D \t f_{10}^{-1} - Id \| = \| (A - B)^{-1} \cdot (A - B - Id) \| > \frac1{10},
$$
see~\eqref{eq:185}.
\end{proof}

Now we extend maps~$\t f_{ij}$ from the cube~$Q$ to the whole sphere~$M$, see Figure~\ref{fig:sphere}.

\begin{Prop}    \label{prop:fij_exist}
There exist maps~$f_{ij} \colon M \to M$, $ij \in \{0, 1\}^2$, with the following properties:

1) $f_{ij}|_Q = \t f_{ij}$ for each $ij \in \{0, 1\}^2$;

2) The maps~$f_{ij}$ are Morse-Smale \diffeos of $M$;

3) For any $g_{ij}$ close enough to $f_{ij}$ for almost every $(\o, x) \in X$ there exists~$k \in \bbN$ \st for each $m > k$
$$
\pi G^m (\o, x) \in Q^+.
$$
\end{Prop}

\begin{proof}
As each of the \diffeos $\t f_{ij}$ sends $Q$ strictly inside itself, we can pick the \diffeos $f_{ij}$ \st

a) $f_{ij}|_Q = \t f_{ij}$ for each $ij \in \{0, 1\}^2$;

b) \tes a closed ball $J \subset M \setminus Q$ \st each of $f_{ij}$ is uniformly expanding on $J$;

c) each of $f_{ij}$ has a unique fixed point~$p_{ij}$ outside of $Q$, moreover, $p_{ij} \in J$;

d) \tes $N \in \bbN$ \st $\forall x \in M \setminus (Q \cup J)$ $\forall \om$ $\pi F^N (\om, x) \in Q$.

We want to emphasize here that the construction of the whole map~$F \colon X \to X$ is now complete.

\begin{figure}[hbt]
    \centering
    \includegraphics[width=300pt]{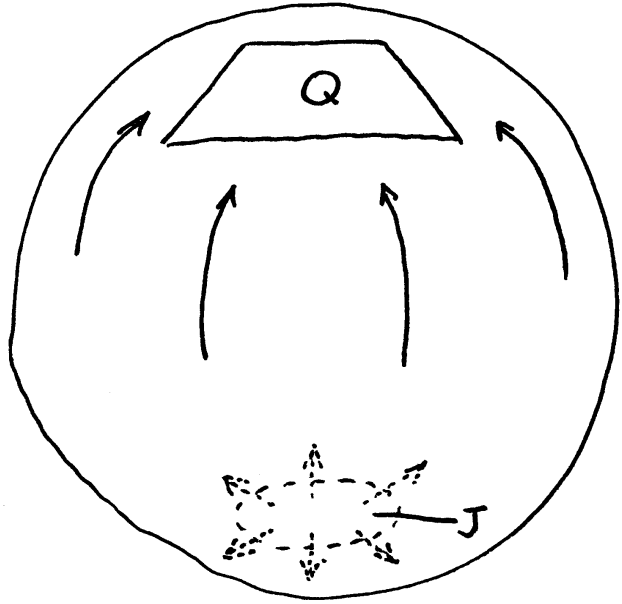}
    \caption{Dynamics of $f_{ij}$ on the sphere~$M$}
    \label{fig:sphere}
\end{figure}

The properties~1 and~2 immediately follow from these conditions. In order to prove property~3 we employ the idea of Lemma~1 from~\cite{bib:IN}. Consider first the inverse map~$F^{-1}$ restricted to~$X^- := \Stt \times J$:
$$
F^{-1}: (\o, x) \mapsto (\s^{-1} \o, f^{-1}_{\o_{-1}} (x)).
$$
All the fiber maps are contracting on $J$. Consider the maximal attractor of~$F^{-1}|_{X^-}$:
$$
S := \bigcap_{n > 0} F^{-n} (X^-).
$$
It is a repelling set for $F$. Let
$$
S_m := \bigcap_{k = 0}^m F^{-k} (X^-).
$$
For any $\o \in \Stt$ we denote
$$
S_{m, \o} := \{ x \in J \,|\, (\o, x) \in S_m \}.
$$
Let $l < 1$ be the contraction coefficient of all the fiber maps on~$J$. Then
$$
\diam S_{n, \o} \le l^n \diam J.
$$
Hence, the intersection of all the nested compact sets $S_{m, \o}$ is one point. Denote it by $\gamma (\o)$.

Thus $S$ is the graph of a function~$\gamma \colon \Stt \to J$. It intersects each fiber exactly at one point. By Fubini Theorem, its measure equals zero.

The maximal attractor $S$ of $F^{-1}|_{X^-}$ consists of all complete orbits of this map. Any other point has a finite past orbit under~$F^{-1}$. This implies that for any point~$p$ from $X^- \setminus S$ there exists $k \in \bbN$ such that $F^k(p) \notin X^-$. As $\mu_X S = 0$, this proves that almost every orbit leaves~$X^-$.


Now condition d) implies that almost every orbit once gets into $\pi^{-1} Q$ and thus into $\pi^{-1} Q^+$, see \eqref{eq:astat_detail}, \eqref{eq:def_fg}, \eqref{eq:def_fij}, \eqref{eq:def_f10}. Due to the same equations the orbit never comes out of $\pi^{-1} Q^+$, so property~3 is proved for $f_{ij}$. But the conditions b)~---d) are open. Thus for any $g_{ij}$ close enough to $f_{ij}$ we also have property~3.
\end{proof}

\begin{Rem} \label{rem:4}
Proposition~\ref{prop:fij_exist} implies $\pi \Ast(G) \subset Q^+$.
\end{Rem} 
\section{The invisibility of the set $R$} \label{sec:4}

First we prove Theorem~\ref{thm:1} for step \skpr $F$. This section deals with the invisibility of $R$ and the next one establishes property~\eqref{eq:astat_detail} of the statistical attractor.

We show that in order to bring a fiber point into~$R$, one has to meet in the base an extraordinary rare word $\o_1^2 \ldots \o_{n^2}^2$ consisting of $n^2$ consecutive zeros. Thus the invisibility rate~$\e$ for $R$ is not greater than $2^{-n^2}$.

Denote
\begin{equation}    \label{eq:w}
    W' := \{\,(x,y) \in W \,|\, y \le \frac14 - h\} = \left[\frac1n, \frac4n\right] \times \left[-\frac2n, \frac14 - h\right].
\end{equation}

\begin{Lem} \label{lem:4}
Let $k > n$ and $\pi F^k (\o, x) \in W'$.
Then
\begin{equation}    \label{eq:z1}
    ( \o^1_{k-n} \ldots \o^1_{k-1} ) = ( 0 \ldots 0 ),
\end{equation}
\begin{equation}    \label{eq:z2}
    ( \o^2_{k-n} \ldots \o^2_{k-1} ) = ( 0 \ldots 0 ).
\end{equation}
\end{Lem}

\begin{proof}
Here we use the same argument as in~\cite[Proposition 4]{bib:IN}. First we prove the part about $\o^1$. Let $j$ be the position of the last occurrence of $1$ in the sequence $\o^1$ before $k$:
$$
j = \max \{ i < k \,|\, \o^1_i = 1 \}.
$$
If there is no such $j$ then the Lemma is proved because $k > n$. Remember that $\pi_1 \colon X \to S^1$ is the projection onto the fiber's first coordinate:
$$
\pi_1: (\o, x^1, x^2) \mapsto x^1.
$$
Since $\o^1_j = 1$, we have
$$
x^1_j := \pi_1 F^{j+1} (\o, x) = f_1 ( \pi_1 F^{j} (\o, x) ) > f_1 (0) = \frac13.
$$
Due to the choice of the map~$f_0$, for all $l = 1, \ldots, k-j$
$$
x^1_{j+l} := \pi_1 F^{j+1+l} (\o, x) > x^1_j - l \cdot \frac{1}{8n} > \frac13 - l \cdot \frac{1}{8n}.
$$
But $\pi F^k (\o, x) \in W$ means $x^1_k \le \frac18$. Thus
$$
\frac13 - (k-j) \cdot \frac1{8n} \le \frac18,
$$
which implies
$$
k-j > n.
$$

Now we established that the last $n$ symbols in $\o^1$ are zeros. This means that the last $n$ fiber maps which brought a point $x$ into $W'$ were either $f_{00}$ or $f_{01}$. But the map $f_{01}$ sends the whole square $[0, 1]^2$ into the upper rectangle $P$, see~\eqref{eq:def_p}, which is invariant under both $f_{00}$ and $f_{01}$ and has empty intersection with the region~$W'$. Thus all the last $n$ fiber maps had to be $f_{00}$'s.
\end{proof}

We denote by $A_m$, $m = 1, \ldots, 4n$, the rectangular regions of height~$h$
$$
A_m = \{(x, y) \in Q^+ \,|\, 1/4 - mh \le y < 1/4 - (m-1)h \}.
$$
The partition of the lower part of $Q^+$ by these regions is shown on Figure~\ref{fig:fij} for $\t f_{00}$ and $\hat f_{10}$.

In the following two Propositions, we study the dynamical behavior of the regions $A_m$ under the maps $f_{00}$ and $f_{10}$. The results are then used in Lemma~\ref{lem:3}.

\begin{Prop}    \label{prop:3.1}
The regions~$A_m$ are invariant under the map~$f_{00}$,
$$
f_{00} (A_m) \subset A_m,
$$
moreover,
$$
f_{10} (A_m \setminus W) \subset A_m.
$$
\end{Prop}

\begin{proof}
The top and bottom sides of these regions are segments of the invariant manifolds of the map~$f_{00}$, hence the regions themselves are (forward) invariant under this map. The same reason works for the restriction of the map~$f_{10}$ to the region~$A_m \setminus W$.
\end{proof}

\begin{Prop}    \label{prop:3.2}
The map~$f_{10}$ in the weak fall-down region $W$ moves points down not more than two regions~$A_m$ at a time:
$$
f_{10} (A_m \cap W) \cap \bigcup_{l > m+2} A_l = \emptyset.
$$
\end{Prop}

\begin{proof}
As $\|\a \b\|_{C^0} = \frac1{10n} < 2h$, this statement follows from the definition of~$f_{10}$, see~\eqref{eq:def_f10} and Proposition~\ref{prop:ab_exist}.
\end{proof}

\begin{Lem} \label{lem:3}
Let $k > n^2$ and $F^k (\o, x) \in R$. Then
\begin{equation}    \label{eq:z2sq}
( \o^2_{k-n^2} \ldots \o^2_{k-1} ) = ( 0 \ldots 0 ).
\end{equation}
\end{Lem}

\begin{proof}
Let $j$ be the position of the last occurrence of $1$ in the sequence $\o^2$ before $k$:
$$
j = \max \{ i < k \,|\, \o^2_i = 1 \}.
$$
If there is no such $j$ then the Lemma is proved, because $k > n^2$. Let $j < k_1 < \ldots < k_m \le k-1$ be the positions such that $\o^1_{k_l} \o^2_{k_l} = 10$ and for $(\t \o, \t x) = F^{k_l} (\o, x)$ we have $\t x \in W$. As $\o^2_j = 1$, the last fiber map in $F^{j+1}$ is either $f_{01}$ or $f_{11}$, so
$$
\pi F^{j+1} (\o, x) \in P.
$$
Propositions~\ref{prop:3.1} and \ref{prop:3.2} imply $m \ge \frac{2n}{2} = n$. Lemma~\ref{lem:4} gives us $k_l - k_{l-1} > n$ $\forall l = 1, \ldots, m$. Summarizing these statements we obtain
$$
k-1-j > \sum_{l = 1}^{m} (k_l - k_{l-1}) > n^2.
$$
Hence $\o^2_j = 0$ for all $k-n^2 \le i \le k-1$.
\end{proof}

Now we are ready to complete the proof of Theorem~\ref{thm:1} for the single map~$F$, by proving that the set~$R$ is $\e$-invisible. Almost every point~$(\o, x)$ visits~$R$ with the frequency not greater than the occurrence of $n^2$ consecutive zeros in the sequence~$\o^2$. By the ergodicity of the Bernoulli shift, for almost all~$\o$, this frequency equals
$$
\e = 2^{-n^2}.
$$ 
\section{The statistical attractor} \label{sec:5}

In this section we prove the first part of Theorem~\ref{thm:1}:
\begin{Lem}     \label{l:main}
For the skew product~$F$ defined above,
\begin{equation}    \label{e:stat}
Q^- \subset \pi \Ast(F).
\end{equation}
\end{Lem}
The right inclusion in~\eqref{eq:astat_detail} is already justified by Remark~\ref{rem:4} to Proposition~\ref{prop:fij_exist}:
$$
\pi \Ast(F) \subset Q^+.
$$
In the following two subsections we establish several lemmas which are key tools for the study of the statistical attractor.

\subsection{Hutchinson lemma and its modifications} \label{ss:5.1}

Let $\A$ be any finite alphabet and the set of maps $f_{\a}$ be indexed by $\a \in \A$. Let $w$ be any finite word $\o_1 \ldots \o_m$, $\o_l \in \A$. Then we denote
$$
f_w = f_{\o_m} \circ \dots \circ f_{\o_1}.
$$

The following lemma is due to Hutchinson~\cite{bib:Hut}.

\begin{Lem} [Hutchinson]   \label{lem:hut}
Consider a metric space~$(M, \r)$ and maps $f_\a \colon M \to M$, $\a \in \A$. Suppose \te compact sets $K \subset K^+ \subset M$ \st
\begin{description}
\item[(invariance)] $\forall \a \in \A, \quad f_\a (K^+) \subset K^+$;
\item[(contraction)] $\exists \, 0 < \l < 1 \ \text{ such that } \ \forall \a \in \A, \ \forall x,y \in K^+,$
$$
\r (f_\a(x), f_\a(y)) < \l \r (x, y);
$$
\item[(coverage)] $\bigcup\limits_{\a \in \A} f_\a(K) \supset K$.
\end{description}
Then for any open $U$, $U \cap K \ne \emptyset$, \tes a finite word
$$
w = \o_1 \ldots \o_m, \ \o_i \in \A
$$
such that the corresponding composition of maps brings the whole $K^+$ into $U$:
$$
f_w (K^+) \subset U.
$$
\end{Lem}

We call the word~$w$ \emph{a critical word for $U$}.

The idea of the proof for Hutchinson lemma is so transparent that we decided to include the proof here.

\begin{proof}
Fix any $x \in K$. By the coverage assumption we can choose $\o_1 \in \A$ such that $x \in f_{\o_1} (K)$ which is equivalent to $f_{\o_1}^{-1} (x) \in K$. Then we can choose $\o_2 \in \A$ \st $x \in f_{\o_1} \circ f_{\o_2} (K)$. By induction, we obtain a sequence $(\o_n)$ \st $\forall m \in \bbN$ and $w_m = \o_m \ldots \o_1$ (note that $\o_m$ goes \emph{first} here) we have $x \in f_{w_m} (K)$. This implies $x \in f_{w_m} (K^+) = K_m$.
Now remember that each of the maps~$f_\a$ is uniformly contracting on $K^+$ and the set~$K^+$ is invariant under these maps. Thus, the diameters of sets~$K_m$ tend to zero and \tes $m \in \bbN$ \st $K_m \subset U$. The word $w_m$ is the word we looked for.
\end{proof}

Now we develop two modifications of Hutchinson lemma.

\begin{Lem}[Robust Hutchinson lemma]   \label{lem:robust}
Consider a Riemannian manifold~$(M,\r)$
and homeomorphisms~$f_\a \colon M \to f_\a(M)$.
Suppose \te compact sets $K^- \subset K \subset K^+ \subset M$ \st
\begin{description}
\item[(robust inclusion)] $K^-$ belongs to $K$ together with $\e$-neighborhood~$U_\e(K^-)$;
\item[(robust invariance)] $\forall \a \in \A, \quad f_\a (K^+) \subset int(K^+)$;
\item[(contraction)] $\exists \, 0 < \l < 1 \ \text{ such that } \ \forall \a \in \A, \ \forall x,y \in K^+,$
$$
\r (f_\a(x), f_\a(y)) < \l \r (x, y);
$$
\item[(robust coverage)] $\bigcup\limits_{\a\in\A} f_\a (K^-) \supset K$.
\end{description}
Then \tes $\delta > 0$ \st for any set of maps~$\{g_\a \,|\, \a \in \A\}$ that are
\begin{itemize}
\item $\delta$-close to~$f_\a$ in $C^0(M)$,
\item contracting in $K^+$,
\end{itemize}
and for $K \subset K^+$ the assumptions of Hutchinson lemma hold.

Moreover, if $\mlip f_\a^{-1} < L$, then $\delta > \frac{\e}{2L}$.
\end{Lem}

\begin{proof}
The robust invariance property survives small $C^0$-perturbation.
By assumption, the maps~$g_\a$ are contracting on~$K^+$. It remains to prove that
$$
\bigcup\limits_{\a\in\A} g_\a(K) \supset K.
$$
For this it is sufficient to prove that
\begin{equation}    \label{e:incl}
\forall \a\in\A \quad f_\a(K^-) \subset g_\a(K)
\end{equation}
Robust inclusion assumption, compactness of $f_\a(K^-)$ and equality $\partial f_\a(K) = f_\a (\partial K)$ for any $\a \in A$ imply that \tes $\e' > 0$ \st
$$
U_{\e'} (f_\a(K^-)) \subset f_\a (K).
$$
Thus for any $g_\a$ that are $\frac{\e'}{2}$-close to $f_\a$ in $C^0(M)$, we have
$$
U_{\e' /2}(f_\a(K^-)) \subset g_\a(K),
$$
which implies~\eqref{e:incl}. Finally, if $\mlip f_\a^{-1} < L$, then
$$
\forall x,y \quad \frac{1}{L} \r(x,y) < \r(f_\a(x),f_\a(y)).
$$
This implies that $\e'$ and $\delta$ can be taken so that $\e' > \frac{\e}{L}$, $\delta > \frac{\e}{2L}$.
\end{proof}

\begin{Lem} [Robust Hutchinson lemma for Cartesian products]   \label{lem:cart}
Consider two Riemannian manifolds~$(M, \rho)$ and $(N, d)$, and two sets of homeomorphisms $\{f_\alpha \colon M \to M \,|\, \a \in \A \}$ and $\{g_\beta \colon N \to N \,|\, \beta \in  \mathcal B \}$.
Suppose the maps $g_\beta$ satisfy assumptions of Robust Hutchinson lemma for the sets  $L^- \subset L \subset L^+ \subset N$. Suppose there exists a collection of words $w_j, j = 1, \ldots, m$ in the alphabet $\mathcal A$ such that the maps $f_{w_j}$ satisfy assumptions of Robust Hutchinson lemma for the sets $K^- \subset K \subset K^+ \subset M$.

Let $B_j$ be the set of all words in the alphabet $\mathcal B $ of length $|w_j|$ (length of $w_j$). Then the maps $\{f_{w_j v} := f_{w_j} \times g_{v} \,|\, j \in \{1, \ldots, m\}, \ v \in B_j\}$ satisfy the assumptions of Robust Hutchinson lemma for the domains $\mathcal K^- = K^- \times L^-$, $\mathcal K = K \times L$ and $\mathcal K^+ = K^+ \times L^+$.
\end{Lem}

\begin{proof}
The domain~$\K^+$ is robustly invariant for~$f_{w_j v}$, because
$$
\forall j \ f_{w_j}(K^+) \subset int K^+, \quad \forall\beta \ g_\beta(L^+) \subset int L^+.
$$
The maps~$f_{w_j v}$ are contracting on~$\K^+$ because the maps~$f_{w_j}$ and $g_\beta$ are contracting on~$K^+$ and~$L^+$ respectively.

The images of $\K^-$ under~$f_{w_j v}$ cover~$\K$. Indeed, for any $j$ we have
$$
\bigcup\limits_{v \in B_j} g_v(L^-) \supset L,
$$
which implies
$$
\bigcup\limits_{j = 1}^{m} \bigcup\limits_{v\in\B_j} f_{w_j v} (K^- \times L^-) = \bigcup\limits_{j=1}^{m} \left( f_{w_j}(K^-) \times \bigcup\limits_{v\in\B_j} g_v(L^-) \right)\supset \bigcup\limits_{j=1}^{m} f_{w_j}(K^-) \times L \supset \K.
$$
\end{proof}

\subsection{Critical words} \label{ss:5.2}

\begin{Lem} [A. Negut]   \label{lem:negut}
Consider a skew product $F$ over a classical Bernoulli shift $\s$ with a fiber $M$:
$$
F: (\o, x) \mapsto (\s\o, f_{\o_0}(x)).
$$
Let $\pi$ be the natural projection
$$
\pi \colon X = \S^2 \times M \to M.
$$
Let $K \subset M$ be an open subset, for which
$$
\pi \Ast (F) \subset K.
$$
Moreover, $K$ is an absorbing set for each~$f_j$:
$$
\forall j \ f_j (K) \subset K.
$$
Suppose that for any open \nbd $U$ of a point $x \in M$ \tes a critical word $w = \o_1 \ldots \o_m$ such that
$$
f_w(K) \subset U.
$$
Then $x \in \pi \Ast (F)$.
\end{Lem}

\begin{proof}
Let us prove that for any \nbd $U$ of $x$, the set~$\pi^{-1} U$ is visited by almost all points with positive frequency. By Definition~\ref{def:astat}, this implies that $x \in \pi \Ast (F)$.

Take any point $(\o, y) \in \pi^{-1} (K)$. Suppose that $\o$ contains the critical word of length $m$ at position $k$, that is:
$$
\o_{k} \dots \o_{k+m-1} = w.
$$
Let
$$
(\s^k \o, y') = F^k (\o, y).
$$
Then $y' \in K$, because $y \in K$ and $\forall j \ f_j (K) \subset K$. Now,
$$
F^{k+m} (\o, y) = (\s^{k+m} \o, f_w (y')) \in \pi^{-1} U
$$
by the choice of word~$w$. Hence, any occurrence of a subword $w$ in $\o$ corresponds to a visit of a point $(\o, y)$ to $\pi^{-1} U$. By the ergodicity of Bernoulli shift, any word of length $m$ is met in a typical sequence $\o$ with an average frequency of $2^{-m}$. Hence, almost all points from $\pi^{-1} K$ visit $\pi^{-1} U$ with a positive frequency. On the other hand, as $K$ is a \nbd of $\pi \Ast (F)$, almost all points of $X$ visit $\pi^{-1} K$ with a positive frequency. Hence,
$$
\Ast (F) \cap \pi^{-1} U \neq \emptyset.
$$
This implies the Lemma.
\end{proof}

\subsection{Finding the critical words} \label{ss:finding}

In this subsection we prove Lemma~\ref{l:main} using the techniques developed in Subsections~\ref{ss:5.1}~--\ref{ss:5.2}.

\begin{Lem} \label{lem:inclu}
For any $x \in Q^- \subset Q$ and any \nbd $U$ of $x$, \tes a critical word $w$ such that
\begin{equation}    \label{eq:subsetu}
f_w (Q^+) \subset U.
\end{equation}
\end{Lem}

Together with Lemma~\ref{lem:negut} and Proposition~\ref{prop:fij_exist} (property~3) this implies~\eqref{e:stat}.

\begin{proof}
We are going to employ the Robust Hutchinson lemma for Cartesian products to obtain the critical words we need. Unfortunately, the maps~$f_{ij}$ are not Cartesian products on the whole~$Q^+$, neither they are contracting. And the regions where they are contracting Cartesian products, are not invariant under the dynamics.

We will now introduce some new maps which are the combinations of~$f_{ij}$, and construct the sets~$\mathcal K^- \subset \mathcal K \subset \mathcal K^+ \subset Q^+$ with the following properties. On the one hand, these new maps and $\mathcal K^-, \mathcal K, \mathcal K^+$ satisfy the assumptions of Robust Hutchinson lemma for Cartesian products. Thus for any open~$U$ that intersects~$\mathcal K$ \tes a word $w$ \st
\begin{equation}    \label{e:incl1}
f_w (\mathcal K^+) \subset U.
\end{equation}
On the other hand, we will prove that there exists a word $w$ such that $f_w (Q^+) \subset \mathcal K^+.$
At last, we will prove that for any $z \in Q^-$ \te $y \in \mathcal K$ and a word $w''$ \st
$$
f_{w''} (y) = z.
$$
All together, this will imply~\eqref{eq:subsetu} for any $x \in Q^-$ and open set $U \ni x$.

The maps~$f_{ij}$, $ij \in \{0,1\}$, are Cartesian products on $[d,1] \times [0,1]$, $d = \frac5n$:
$$
f_{ij} = f_i \times g_j.
$$
The maps~$g_0$, $g_1$ are contracting on~$\left[\frac14 - \frac{h}{2}, 1 + \frac{h}{2}\right] =: J^+$. The maps~$g_0, g_1$ bring $J^+$ strictly into itself. On the other hand, for $J := \left[\frac14 + \frac{h}{2}, 1 - \frac{h}{2} \right]$, and for $\delta$ small we have:
$$
g_0(J^-) \cup g_1(J^-) \supset J,
$$
where~$J^- := \left[ \frac14+\frac{h}{2}+\delta, 1 - \frac{h}{2}-\delta\right]$. So $g_0, g_1, J^-, J, J^+$ satisfy the assumptions of Lemma~\ref{lem:cart}.

We will now find segments~$L^- \subset L \subset L^+$ and combinations of maps~$f_0, f_1$ that satisfy the assumptions of Lemma~\ref{lem:cart}. Recall that on $I = [-1, 2]$ the maps $f_0, f_1$ are linear contractions with fixed points $0$ and $1$ and coefficients $\lambda = 1 - \frac{1}{8n}$ and $\mu = \frac 2 3$ respectively.

Denote~$f_1(0) = a$. Let us fix any $k$ and $c$ \st
\begin{equation}    \label{e:ck}
f_0^k (a) =: c \in [d, 2d).
\end{equation}
Let:
$$
L^+ = [a, 1 + c], \ L = [a + c, 1 - \frac 3 4 c], \ L^- = [a + 2c, \ 1 - c].
$$
The maps $f_\alpha $ are taken as
\begin{equation}    \label{e:hj}
h_m = f_1 \circ f_0^m, \ m = 0, \dots , k.
\end{equation}
Now we verify the assumptions of Lemma~\ref{lem:cart}. We obviously have the robust inclusion property.

As $\forall m \ f_0^m(a) > 0$, we have $h_m(a) > a$. Also $\forall x \in (1,1+c]$ we have $f_0(x) < x$, $f_1(x) < x$. Hence, $h_m(1+c) < 1+c$. This proves the robust invariance of~$L^+$.

The maps $h_m$ are contracting on $L^+$, because $f_0, f_1$ are contracting on $I^+$.

Note that~$\lambda^k = \frac{c}{a}$, see~\eqref{eq:def_fg} and \eqref{e:ck}. Let~$\bigcup_m h_m(L^-) = f_1(L')$. We have:
$$
L' = \bigcup_{m=0}^k f_0^m(L) = \left[c+\frac{2c^2}{a}, 1-c\right].
$$
Then
$$
f_1(L') = \left[ a + \frac23 \left( c + \frac{2c^2}{a} \right), 1 - \frac23 c \right] \supset L.
$$
This gives us the robust coverage.

Therefore we can apply Lemma~\ref{lem:cart} (Robust Hutchinson lemma
for Cartesian products), to
$$
\mathcal K^- = L^- \times J^-, \ \mathcal K = L \times J, \ \mathcal K^+ = L^+ \times J^+,
$$
and the maps
\begin{equation}    \label{e:f_jv}
\{ f_{m,v} := h_m \times g_{v}\,|\, |v| = m \}, \ m = 0, \dots , k.
\end{equation}
These maps are compositions of the maps~$f_{ij}$, because $\mathcal K^+ \cap W = \emptyset$. Thus we obtain~\eqref{e:incl1}.

\begin{Rem} \label{r:perturb}
We also have~\eqref{e:incl1} for perturbations of~$f_{ij}$, \st the perturbations of all the compositions~\eqref{e:hj} are small enough in~$C^0$ simultaneously. In Section~\ref{sec:8} we will show it is possible.
\end{Rem}

\begin{Prop}    \label{prop:crit_word_for_K}
There exists a word $w$ such that
\begin{equation}    \label{eq:inclu1}
f_w (Q^+) \subset \mathcal K^+.
\end{equation}
\end{Prop}

\begin{proof}
The following argument is illustrated by Figure~\ref{fig:inclu_fix}. Consider $y = f_{00}(1, 1)$. By definition of the map~$f_{00}$, see~\eqref{eq:def_fij} and Figure~\ref{fig:fij}, the point $y$ together with some \nbd $U_0$ is within $\mathcal K$.

On the other hand, $f_{11}$ is contracting on $Q^+$ with the unique attractor~$(1, 1)$. Let us take $m$ so large that $f_{00} \circ f_{11}^m (Q^+) \subset U_0$, see Figure~\ref{fig:inclu_fix}. The word
$$
w := \underbrace{(11) \ldots (11)}_{m} (00)
$$
is the desired one.
\end{proof}

\begin{figure}[hbt]
    \centering
    \includegraphics[width=300pt]{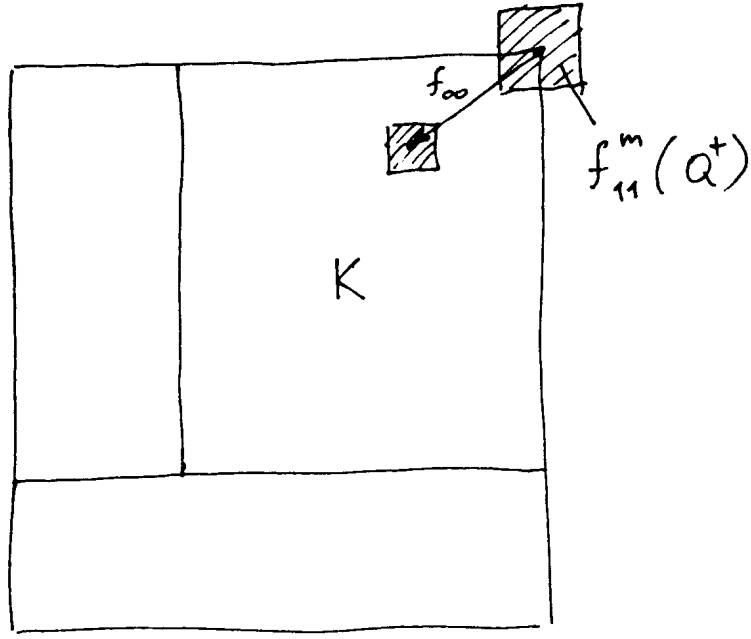}
    \caption{Getting inside of~$K$}
    \label{fig:inclu_fix}
\end{figure}

\begin{Prop}    \label{prop:7}
For every $x \in Q^- \cap P$ \tes a word $w$ \st $f_{w}^{-1}(x) \in \mathcal K$.
\end{Prop}

\begin{proof}
Consider the sequence
$$
\mathcal P_n = Q^- \cap \bigcup_{ij} f_{ij} (\mathcal P_{n-1}), \quad \mathcal P_0 = K.
$$
It is a monotonous sequence of rectangles that tends to~$Q^- \cap P$. Thus the word~$w$ exists.
\end{proof}

\begin{Prop}    \label{prop:9}
For any $m = 1, \ldots, 4n$ for every~$q \in A_m$ \tes a word~$w(q)$ such that~$f_{w(q)}^{-1} (x) \in Q^- \cap P$.
\end{Prop}

\begin{proof}
The construction of the word~$w(q)$ for $q \in A_m$, see Figure~\ref{fig:fij}, is done by induction in $m$. Assume that such a word exists for any $q \in A_l$, $l < m$. Consider the region~$B = f_{10} (D)$, see~\eqref{eq:def_dw}. The equations~\eqref{eq:def_fg} and \eqref{eq:def_f10} imply that $B$ is a rectangle of width~$\frac2{3n}$, its sides parallel to horizontal and vertical axes. Note that for any $m = 1, \ldots, 4n$ the set $A_m \setminus B$ has two connected components:
$$
A_m \setminus B = L_m \sqcup R_m,
$$
where $L_m$ stands for the left one and $R_m$ for the right one. For a point $q \in A_m$ we consider $3$ cases:

1) $q \in A_m \cap B$. Then $\pi_2 (f_{10}^{-1} q) - \pi_2 q > h$, see~\eqref{eq:projk}, \eqref{eq:def_f10}, \eqref{eq:def_dw}, which implies that $ f_{10}^{-1} q$ lies in~$D$ above the region $A_m$. Additionally, if $f_{10}^{-1} q \in P$, \tes a word~$w'$, composed of $f_{00}$ and $f_{01}$, \st $f_{w'}^{-1} \circ f_{10}^{-1} q \in Q^- \cap P$.

2) $q \in L_m$. Consider the backward orbit of $q$ under the map $f_{00}$:
$$
q_k := f_{00}^{-k} (q), \quad k \in \bbN.
$$
The equations~\eqref{eq:def_fg} and~\eqref{eq:def_fij} imply that
$$
0 < \d(q) < \pi_1 q_k - \pi_1 q_{k-1} < \frac1{8n},
$$
where $\d(q)$ depends only on the initial point $q$. As the width of the stripe $B$ is greater than~$\frac1{8n}$, \tes a $k \in \bbN$ \st $ f_{00}^{-k} q \in B$ and we are in the settings of case~1.

3) $q \in R_m$. It follows from~\eqref{eq:def_fg} and \eqref{eq:def_f10} that for each $(x, y) \in R_m$ we have $x > \frac13$. Let $q'_k := f_{10}^{-k} (q)$. Consider their abscissas $\pi_1 q'_k$ on a logarithmic scale centered at $x = 1$. The equations~\eqref{eq:def_fg} and~\eqref{eq:def_f10} imply that for any $k \in \bbN$
$$
\log (1 - \pi_1 q'_k) - \log (1 - \pi_1 q'_{k-1}) = \log \frac32,
$$
while
$$
\log (1 - 0) - \log (1 - \frac23) > \log \frac32.
$$
Also note that for any $q' \in R_m$ we have $f_{10}^{-1} q' \in A_m$. Thus \tes $k \in \bbN$ \st $q'_k \in A_m \setminus R_m$. So we have just reduced the case~3 to the cases~1 and~2. This proves the Proposition.
\end{proof}

This Proposition completes the proof of Lemma~\ref{l:main}. Thus we have proved the first statement of Theorem~\ref{thm:1} for the single map $F$.
\end{proof}

\section{Perturbation in the space of step skew products}   \label{sec:6}

Now we are going to prove Theorem~\ref{thm:1} for any step skew product~$G$ that is close enough to $F$. The distance in the space of step \skprs is always interpreted as~\eqref{eq:def_metric}.

In this section, we establish some basic facts about the dynamics of $g_{ij}$ for any~$G$ which is close to the initial map~$F$.

\begin{Prop}    \label{prop:struct}
For any~$G$ close enough to~$F$ for every $ij \in \{0, 1\}^2$ \tes \homeo $H_{ij} \colon M \to M$ \st the following diagram commutes:
$$
\begin{CD}
M @>{f_{ij}}>> M \\
@V{H_{ij}}VV            @VV{H_{ij}}V \\
M @>{g_{ij}}>> M
\end{CD}
$$
and
\begin{equation}    \label{eq:r}
d_{C^0}(H_{ij}, Id) \le r, \quad d_{C^0}(H_{ij}^{-1},Id) \le r,
\end{equation}
where $r = \frac{h \nu}{10} = \frac{1}{160n^2}$, see~\eqref{eq:def_h}.
\end{Prop}

\begin{Rem} Inequalities~\eqref{eq:r} imply $d_{C^0} (f_{ij}, g_{ij}) \le r$, $d_{C^0} (f_{ij}^{-1}, g_{ij}^{-1}) \le r$.
\end{Rem}

\begin{proof}
According to~\eqref{eq:def_metric}, the closeness of $G$ to $F$ is equivalent to the $C^1$-closeness of each pair of fiber maps $g_{ij}^{\pm 1}$ to the corresponding maps~$f_{ij}^{\pm 1}$. Now the claim of the Proposition follows from the fact that the maps $f_{ij}$ are Morse-Smale diffeomorphisms, see section~\ref{sec:2}. These maps are structurally stable. The estimates~\eqref{eq:r} are straightforward and we skip them.
\end{proof}

In the following text we will always assume that the maps $F,G$ are at least that close that Proposition~\ref{prop:struct} works.

We will say that a map~$f \colon M \to M$ \emph{moves the points to the right (to the left)} in some region $E \subset Q$ if \tes $\a > 0$ \st $\forall q \in E$, $q = (x, y)$, $f (x, y) = (x', y')$
$$
x' - x > \a, \ (x - x' > \a).
$$

We will say that a map~$f \colon M \to M$ \emph{moves the points up (down)} in some region $E \subset Q$ if \tes $\a > 0$ \st $\forall q \in E$, $q = (x, y)$, $f (x, y) = (x', y')$
$$
y' - y > \a, \ (y - y' > \a).
$$

\begin{Rem}  \label{rem:push}
The set of \diffeos which move the points in some direction (right, left, up or down) in some compact region is $C^0$-open. Thus if a fiber map~$f_{ij}$ moves the points in some direction, then we have the same for the corresponding fiber map~$g_{ij}$ for any $G$ which is close enough to $F$.
\end{Rem}


Now denote
\begin{equation}    \label{e:rho}
\rho = \frac{r}{\nu} = \frac{h}{10}.
\end{equation}
Note that for any $n > 10$ we have $10r < \rho$. For every~$ij \in \{0, 1\}^2$ let us subtract the $\rho$-\nbds of the invariant manifolds of~$f_{ij}$ from~$Q^+$; for $ij = 10$ we also subtract the region~$W$, see Figure~\ref{fig:fall}. Denote the result by~$\t Q^+_{ij}$. Each set~$\t Q^+_{ij}$ consists of finitely many linearly connected components.

\begin{Prop}    \label{p:stable_movement}
Fix~$ij \in \{0,1\}^2$. Let~$C$ be a connected component of $\t Q^+_{ij}$. Then either $g_{ij}$ moves all the points of~$C$ to the left or it moves all the points of~$C$ to the right; one has the same alternative for up and down directions. These directions coincide with the ones of $f_{ij}$ in the same regions.
\end{Prop}

\begin{proof}
This statement follows directly from~\eqref{eq:def_fg}, \eqref{eq:def_fij}, \eqref{eq:def_f10}, Proposition~\ref{prop:struct} and Remark~\ref{rem:push}.
\end{proof}

Now let us introduce the following notations for certain subsets of $Q^+$. We let
$$
S_m^- := \{ (x, y) \in Q^+ \,|\, d (y, \frac14 - (2m - 2)h ) \le \rho \},
$$
$$
S_m^+ := \{ (x, y) \in Q^+ \,|\, d (y, \frac14 - (2m - 1)h ) \le \rho \},
$$
be the \nbds of the invariant manifolds of $f_{00}$, $S_m^-$ correspond to the strong stable manifolds of the attracting fixed points and $S_m^+$ correspond to the stable manifolds of the saddle points. We also denote
$$
S_{all} := \bigcup\limits_m S_m^+ \cup S_m^-, \quad
U_m := A_{2m-1} \setminus S_{all}, \quad
D_m := A_{2m-2} \setminus S_{all}
$$
Note that~$U_m$, $D_m$ are connected components of $Q^+ \setminus S_{all}$. Here $U$ stands for up, $D$ stands for down, --- the general direction of dynamics by~$f_{00}$ and $f_{10}$ in those regions.

\section{The invisibility in the perturbed skew product} \label{sec:7}

In this section, we prove that the set~$R$, see~\eqref{eq:invis_detail}, is $\e$-invisible, $\e = 2^{-n^2}$, for any~$G$ close enough to~$F$. We follow the strategy of section~\ref{sec:4} where we prove the same property for the single map~$F$.

\begin{Lem} \label{lem:7a}
For any~$G$ close enough to~$F$ we have the following. Let $k > n$ and
\begin{equation}    \label{eq:lem7a}
\pi G^k (\o, x) \in W',
\end{equation}
see~\eqref{eq:w}. Then~\eqref{eq:z1} and~\eqref{eq:z2} hold.
\end{Lem}

\begin{proof}
Let the distance $d(F, G)$ be so small that $d_{C^0} (f_{ij}, g_{ij}) < \frac1{32n}$ for all $ij \in \{0, 1\}^2$. Then for
$$
(x_{ij}, y_{ij}) := g_{ij} (x, y)
$$
the following inequalities hold
$$
x - x_{0j} < \frac1{8n} + \frac1{32n} < \frac1{6n},
$$
$$
x_{1j} > \frac13 - \frac1{32n}
$$
for $j \in \{0, 1\}$. In the same way as in Lemma~\ref{lem:4} they imply~\eqref{eq:z1}.

Now let
$$
Y := \{\, (x, y) \in Q^+ \,|\, y \ge \frac13 - \frac1{32n} \}
$$
and note that $g_{01} (Q^+) \subset Y$. On the other hand, for any $n \ge 0$
$$
g_{00} (Y) \cap W' = \emptyset.
$$
Thus if at least one of $\o^2_{k-n}, \ldots, \o^2_{k-1}$ was not zero, \eqref{eq:lem7a} could not be true, therefore we have~\eqref{eq:z2}.
\end{proof}

\begin{figure}[hbt]
    \centering
    \includegraphics[width=300pt]{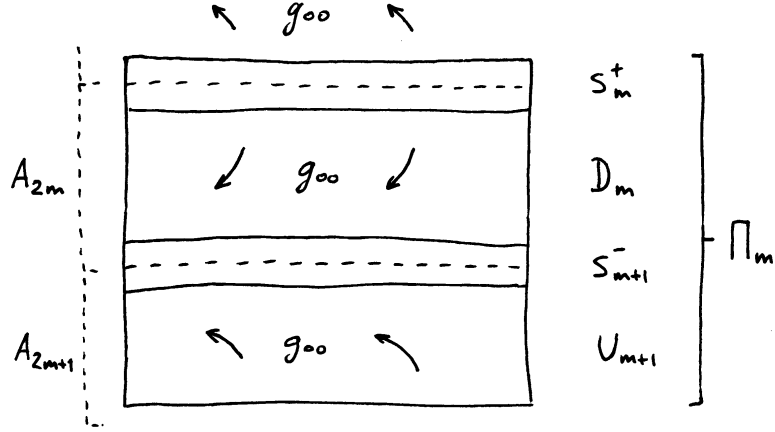}
    \caption{The block~$\Pi_m$ and the map~$g_{00}$}
    \label{fig:pi}
\end{figure}

Now let us break the lower part of $Q^+$ into the following blocks~$\Pi_m$:
$$
\Pi_m := S_{m}^+ \cup D_{m} \cup S_{m+1}^- \cup U_{m+1},
$$
see Figure~\ref{fig:pi}. Each $\Pi_m$ contains two \nbds of invariant manifolds of $f_{00}$ and two gaps~$U_m$ and $D_m$ between them. The block $\Pi_0$ consists of $S^-_1$ and $U_1$ only.

\begin{Prop}    \label{prop:16}
The blocks $\Pi_m$ do not ``go down'' under the maps~$g_{ij}$, $ij \neq 10$; namely:
$$
g_{ij} (\Pi_m) \cap \bigcup_{l > m} \Pi_l = \emptyset, \text{ for } ij \neq 10.
$$
\end{Prop}

\begin{proof}
The images of~$\Pi_m$ under the maps~$g_{01}$, $g_{11}$ have empty intersection with any of $\Pi_m$. Now we study where~$\Pi_m$ goes under $g_{00}$. Note that $S_m^+ \cup D_m \subset A_{2m-1} \cup A_{2m}$ and recall (see Proposition~\ref{prop:3.1}) that the regions $A_m$ are invariant under $ f_{00}$. Thus due to inequality
$$
d_{C^0} (f_{00}, g_{00}) < r < \rho
$$
we have
$$
g_{00} (S_m^+ \cup D_m) \subset U_\rho (A_{2m-1} \cup A_m)
$$
and
$$
U_\rho (A_{2m-1} \cup A_{2m}) \cap \bigcup_{l > m} \Pi_l = \emptyset.
$$
Now, the points of $S_m^-$ move along $y$ axis not farther than by $\rho$ under $f_{00}$ (remember that they are all attracted to the middle line of $S_m^-$ by $f_{00}$ and the height of $S_m^-$ is $2\rho$) and thus not farther than by $\rho + r < h - 2\rho$ under $g_{00}$. Then
$$
g_{00} (S_m^-) \cap \bigcup_{l > m} \Pi_l = \emptyset.
$$
Finally, the points of $U_m$ all move upwards under $f_{00}$ by at least
$$
\frac{h}{3 \pi} \cdot \sin \frac{\pi \rho}{h} > \frac{h}{3 \pi} \cdot \frac{\pi \rho}{2h}  = \frac{\rho}{6} > 2r > r.
$$
Thus $g_{00}$ also moves the points of $U_m$ upwards, which implies
$$
g_{00} (U_m) \cap \bigcup_{k > m} \Pi_k = \emptyset.
$$
To summarize, we obtain the desired relation
$$
g_{00} (\Pi_m) \cap \bigcup_{k > m} \Pi_k = \emptyset.
$$
\end{proof}

\begin{Prop}    
The blocks $\Pi_m$ do not ``go down'' under the restriction of the map~$g_{10}$ to the complement of the region~$W$; namely:
$$
g_{10} (\Pi_m \setminus W) \cap \bigcup_{k > m} \Pi_k = \emptyset.
$$
\end{Prop}

\begin{proof}
The proof is the exact copy of that of Proposition~\ref{prop:16}.
\end{proof}

\begin{Prop}    
The map~$g_{10}$ in the weak fall-down region $W$ can move the points down not lower than by one block:
$$
g_{10} (\Pi_m \cap W) \cap \bigcup_{k > m+1} \Pi_k = \emptyset.
$$
\end{Prop}

\begin{proof}
The height of each block~$\Pi_m$ equals $2h$. The map~$f_{10}$ moves the points down not lower than by $\frac32 h$. Thus the map~$g_{10}$ moves the points down not lower than by $\frac32 h + r < 2h$.
\end{proof}

\begin{Lem} \label{lem:8}
Let $k > n^2$ and $\pi G^k (\o, x) \in R$. Then we have~\eqref{eq:z2sq}.
\end{Lem}

\begin{proof}
The argument is the same to Lemma~\ref{lem:3}.
\end{proof}

Now we are ready to complete the proof of $\e$-invisibility part of Theorem~\ref{thm:1}, by proving that the set~$R$ is $\e$-invisible. Almost every point~$(\o, x)$ visits~$R$ with at most the frequency of occurrence of $n^2$ consecutive zeros in the sequence~$\o^2$. By the ergodicity of the Bernoulli shift, for almost all~$\o$, this frequency equals
$$
\e = 2^{-n^2}.
$$ 
\section{The statistical attractor for the perturbation} \label{sec:8}

Now we turn back to the proof of the left inclusion in statement~\eqref{eq:astat_detail} of Theorem~\ref{thm:1}, namely, that $Q^- \subset \pi \Ast(G)$. We employ mostly the same ideas and techniques as we did in section~\ref{sec:5}. The critical words are explicitly constructed below for small \nbds of every $x \in Q^-$.

\subsection{Upper rectangle $P^-$}

Let $P^- = Q^- \cap \{ y > \frac14+\rho \}$, see~\eqref{e:rho}.
In this subsection we find the critical words for the \nbds of the points of~$P^-$.
We are going to use the strategy and the results of subsection~\ref{ss:finding}, in particular Remark~\ref{r:perturb}.

To apply the methods of Section~\ref{sec:5}, we have to estimate the discrepancy between $g_w$ and $f_w$ for long words~$w$ within the region~$P^-$.

\begin{Lem}    \label{lem:discrepancy}
Let~$p \in P^-$, and $w = a_1 \ldots a_m$ be \st $\forall l \le m$ $f_{a_l} \circ \ldots \circ f_{a_1} (p) \in P^-$. Then
\begin{equation}    \label{eq:new2}
d(f_w (p), g_w (p)) < \rho.
\end{equation}
\end{Lem}
\begin{proof}
Note that in~$P^-$ all the maps~$f_{ij}$ are uniformly contracting with a rate of at least~$\lambda = 1 - \nu$, see Theorem~\ref{thm:1}.
By induction in $l$ we will prove that for every subword~$w_l = a_1 \ldots a_l$ we have
$$
d(f_{w_l} (p), g_{\o_l} (p)) < \lambda^{l-1} r + \ldots + \lambda r + r.
$$
For $l = 1$ we have this due to Proposition~\ref{prop:struct}. Now,
$$
d(f_{w_{l+1}} (p), g_{w_{l+1}} (p)) \le d(f_{w_{l+1}} (p), f_{a_{l+1}} \circ g_{w_{l}} (p)) + d(f_{a_{l+1}} \circ g_{w_{l}} (p), g_{w_{l+1}} (p)) \le
$$
$$
\le \lambda \cdot d(f_{w_{l}} (p), g_{w_{l}} (p)) + r \le \lambda^l r + \ldots + r.
$$
But for any $l$ we have $\lambda^l r + \ldots + r \le r \frac{1}{1-\lambda} = \frac{r}{\nu} = \rho$.
\end{proof}

Recall~$\mathcal K^-, \mathcal K, \mathcal K^+$ from subsection~\ref{ss:finding}.
In~\eqref{e:f_jv} we defined the maps~$f_{j,v}$ which are the compositions of $j$ maps~$f_{ij}$. Denote by~$g_{j,v}$ the same compositions with~$g_{ij}$.

Lemma~\ref{lem:discrepancy} implies that the maps~$g_{j,v}$ are sufficiently close to~$f_{j,v}$ for every~$G$ from Proposition~\ref{prop:struct}. Thus by Lemma~\ref{lem:cart} for any open~$U$ that intersects~$\mathcal K$ \tes a word $w$ \st
\begin{equation}    \label{e:incl1g}
g_w (\mathcal K^+) \subset U.
\end{equation}

\begin{Prop}    \label{prop:8.5}
For any~$G$ close enough to~$F$ there exists a word $w$ such that $g_w (Q^+) \subset \mathcal K^+$
\end{Prop}

\begin{proof}
The construction is almost the same to the one we used in Proposition~\ref{prop:crit_word_for_K}.

Let $d(F,G)$ be so small that we are in the settings of Proposition~\ref{prop:struct}. Then the map $g_{11}$ has an
attracting fixed point $a$ close to $(1,1)$. For any \nbd $U$
of this point \tes $m$ \st $g^m_{11}(Q^+) \subset U$.

Note that for small enough~$r$, we have $g_{00}(a) \in \mathcal K^+$ which implies we can take $U$ so small that $g_{00}(U)
\subset \mathcal K^+$. Now the desired word is $w := \underbrace{(11) \ldots (11)}_{m} (00)$.
\end{proof}

\begin{Prop}    \label{prop:8.7}
For every $q \in P^-$ \tes a word $w$ \st $g_{w}^{-1} q\in \mathcal K$.
\end{Prop}

\begin{proof}
Inside~$P^-$, the maps~$f_{ij}$ are Cartesian products of~$f_i$ and $g_j$.
Thus we can inductively define the rectangles with horizontal and vertical borders
$$
\mathcal K_{n+1} = \bigcup_{ij \in \{0,1\}^2} f_{ij} (\mathcal K_n),
$$
$\mathcal K_0 = \mathcal K$. They are well-defined as long as their iterations stay outside of~$W$. Note that $\forall n$ $\mathcal K_{n+1} \supset \mathcal K_n$.

The same argument we used in Lemma~\ref{lem:discrepancy} gives us that we can also find rectangles
$$
\tilde{\mathcal K}_{n+1} \subset \bigcup_{ij \in \{0,1\}^2} g_{ij} (\tilde{\mathcal K}_n),
$$
such that $\tilde{\mathcal K}_0 = \mathcal K$, $\tilde{\mathcal K}_{n+1} \supset \tilde{\mathcal K}_n$, and the borders of $\mathcal K_n$ and $\tilde{\mathcal K}_n$ differ no more than by~$\rho$.
But
$$
\bigcup_{n \ge 0} \mathcal K_n \supset \left[\frac4n,1\right] \times \left[\frac14, 1\right],
$$
which implies
$$
\bigcup_{n \ge 0} \tilde{\mathcal K}_n \supset P^-.
$$
Thus the desired word exists.
%
%
\end{proof}

\subsection{Lower region}

Now we construct the critical words for~$Q^- \setminus P^-$.

\begin{figure}[hbt]
    \centering
    \includegraphics[width=300pt]{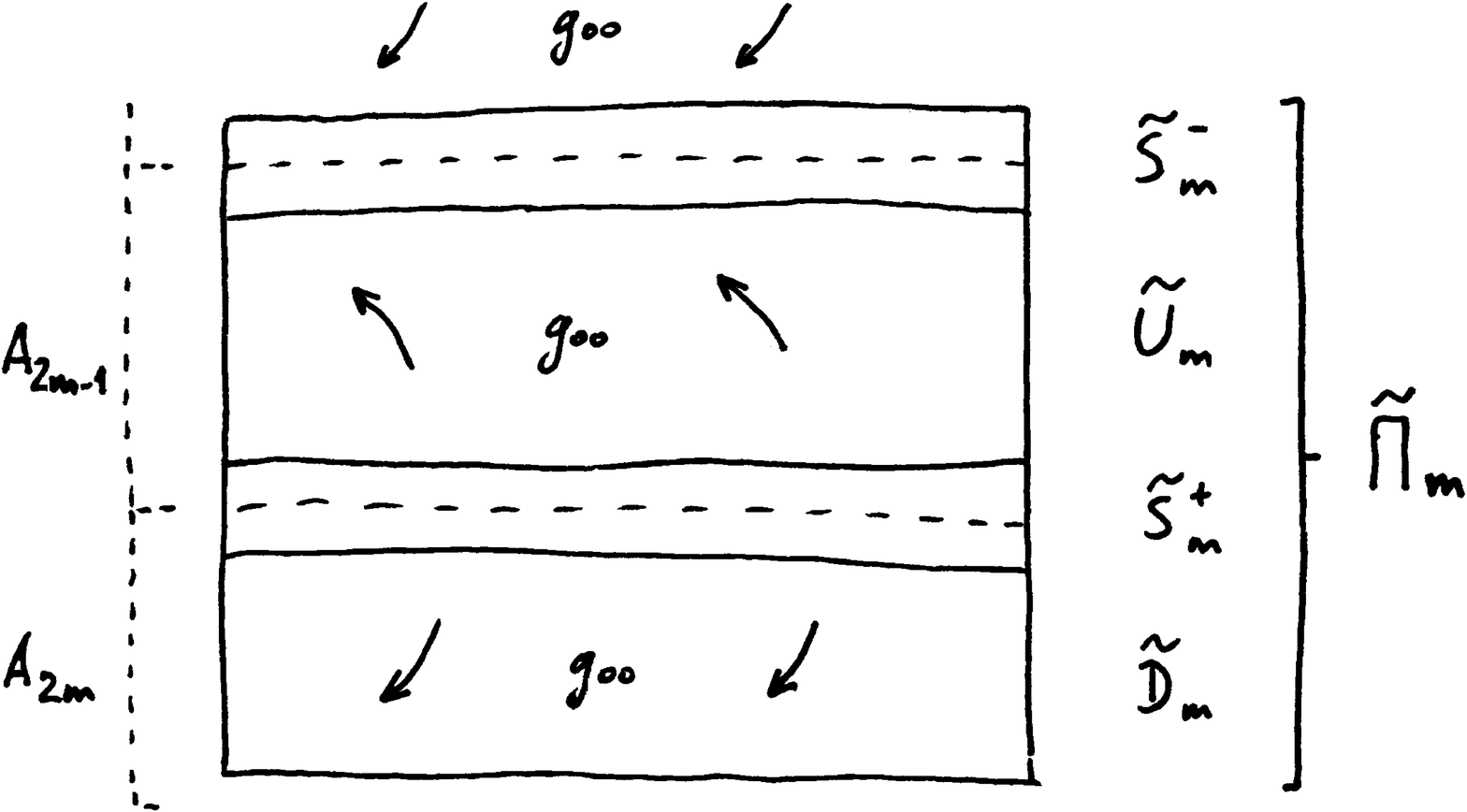}
    \caption{The block $\t \Pi_m$ and the map~$g_{00}$}
    \label{fig:blocks2}
\end{figure}

Let $\t S_m^{\pm} := S_m^{\pm} \cap Q^-$, $\t U_m := U_m \cap Q^-$, $\t D_m := D_m \cap Q^-$, see~\eqref{eq:astat_detail}. Like in section~\ref{sec:7}, it is convenient to break the lower part of $Q^-$ into the following blocks $\t \Pi_m$:
$$
\t \Pi_m := \t S_m^- \cup \t U_m \cup \t S_m^+ \cup \t D_m,
$$
see Figure~\ref{fig:blocks2}. They are similar to the blocks $\Pi_m$ but instead of not going down under forward iterations of $g_{00}$ and $g_{10}$ they do not go down under the backward iterations of these maps, see details below.

\begin{Prop}    \label{prop:8.9}
For any~$G$ close enough to~$F$ we have the following:
$\forall m = 1, \ldots, 2n$ for every~$q \in \t \Pi_m$ and any open \nbd~$U \ni q$ \tes a critical word~$w$ such that
$$
g_w (Q^+) \subset U.
$$
\end{Prop}

\begin{proof}
Like in Proposition~\ref{prop:9}, it is enough to prove that any point~$q \in \t \Pi_m$ can be pushed to the region above $\t \Pi_m$ by the backward iterations of $g_{00}$ and $g_{10}$. We take $d(F,G)$ small enough to be in the settings of Proposition~\ref{prop:struct}.

First of all, note that the regions $Z_m := \t S_m^- \cup \t U_m \cup \t S_m^+$ and $\t D_m$ do not go below themselves under the backward iterations of the maps~$g_{00}$ and $g_{10}$. In fact, according to~\eqref{eq:def_fij}, \eqref{eq:def_f10}, and \eqref{eq:r}, if the backward image of a point of $Z_m$ is below $Z_m$, it can only be within $\t D_m$. But both $g_{00}$ and $g_{10}$ push the points down within $\t D_m$, so this is impossible. Due to the same argument, if the backward image of a point of $\t D_m$ is below $\t D_m$, it can only be within $S^{-}_{m+1}$. But the points of $S^{-}_{m+1}$ either stay within it or go down under $g_{00}$ and $g_{10}$. So this is also impossible.

Consider the region $\t B = g_{10} (D)$. The same consideration as in Proposition~\ref{prop:9} allows us to bring any point $q \in \t \Pi_m$ into $\t \Pi_m \cap \t B$ by the backward iterations of $g_{00}$ and $g_{10}$.

Now, for a point~$q \in \t \Pi_m \cap \t B$ there are two cases.

1) $q \in Z_m$. In this case, the inequality
$$
d_{C^0} (g_{10}^{-1}, f_{10}^{-1}) < r
$$
implies that $g_{10}^{-1} q$ is above $\t \Pi_m$.

2) $q \in \t D_m$. The same inequality tells us that $g_{10}^{-1} q$ lies either within $Z_m$, or is above $\t \Pi_m$.

The proof is completed.
\end{proof}

The assertions of Propositions~\ref{prop:8.7} and \ref{prop:8.9} can be combined into one statement: for any~$G$ close enough to~$F$, we have the following: for every~$x \in Q^-$ and any open \nbd~$U \ni x$, \tes a critical word~$w$ such that
$$
g_w (Q^+) \subset U.
$$

Thus, according to Lemma~\ref{lem:negut} and Remark~\ref{rem:4},
$$
Q^- \subset \pi \Ast(G).
$$
The upper estimate on the projection of the \sa was already given in Remark~\ref{rem:4}:
$$
\pi \Ast(G) \subset Q^+.
$$
This completes the proof of Theorem~\ref{thm:1} about the statistical attractor of~$G$. 
\section{Higher dimension: $k > 2$} \label{sec:9}

In this section we explain how to carry out our construction in the dimension higher than $2$.

\subsection{Construction}

For $k > 2$ we prove Theorem~\ref{thm:1} using induction on $k$. Assume that for a certain $k$ we are able to construct a step \skpr --- the center of the ball~$\mathcal B_p \subset C^1_{k} (L)$ from Theorem~\ref{thm:1}. This step \skpr is uniquely defined by fiber \diffeos $f_{i_1 \dots i_k} \colon S^k \to S^k$, $i_m \in \{0, 1\}$. Our goal now is the construction of fiber \diffeos $f_{i_1 \ldots i_k i_{k+1}}$ for $k+1$-dimensional step \skpr $F$ satisfying Theorem~\ref{thm:1}. They will act on the sphere~$S^{k+1}$.

Within the cube~$Q \subset S^{k+1}$ for each $i_1 \ldots i_k i_{k+1} \neq 1 \ldots 10$ we let
$$
f_{i_1 \dots i_k i_{k+1}} := f_{i_1 \dots i_k} \times g_{i_{k+1}},
$$
which is the direct generalization of~\eqref{eq:def_fij}. Following~\eqref{eq:def_f10}, we also let
$$
f_{1 \ldots 10} := \left( f_{1 \ldots 1} (x_1, \ldots, x_k), g_0(x_{k+1}) - \a(x_k)\b(x_{k+1})\right),
$$
$\a$ and $\b$ are the same as in subsection~\ref{subsec:2-dim}. Now we extend the maps $f_{i_1 \dots i_k i_{k+1}}$ to be the \diffeos of the whole sphere~$S^{k+1}$ like we did in Proposition~\ref{prop:fij_exist}.

\subsection{Proof}

We want to estimate the rate of invisibility of $R = \pi_{k+1}^{-1} \left(-2\nu, \frac1{10}\right)$. The idea of the Sections~\ref{sec:4} and~\ref{sec:7} is employed again: the points of the cube~$Q^+$ come to its upper part quite often and the only chance to go down to~$R$ is to meet an extraordinary rare combination of the letters in base: $n^{k+1}$ consecutive zeros in~$\o^{k+1}$.

Consider the region~$\t R = \pi_{k}^{-1} \left(-2\nu, \frac1{10}\right)$. Note that $\t R$ is $\e$-invisible with $\e \le 2^{-n^k}$ due to the induction hypothesis and to the fact that the first $k$ coordinates of $f_{i_1 \dots i_k i_{k+1}}$ do not depend on $x_{k+1}$.

Now, our $k+1$-dimensional cube~$Q$ is split into $2n+1$ layers~$\Pi_m$, see their analog for~$k = 1$ in section~\ref{sec:7}, which are aligned along the first $k$ coordinates. The transition between the layers in the negative direction (i.e. when $x_{k+1}$ decreases) is possible only when the following two conditions are simultaneously satisfied:

1) $(x_1, \ldots, x_{k+1}) \in \t W \subset \t R$, where $\t W = \{\, x \in Q \,|\, x_k \in \left[ \frac1n, \frac4n \right], x_{k+1} \in \left[ -\frac2n, \frac14 + \frac2n \right] \,\}$;

2) the next fiber map is $f_{1 \ldots 10}$.

Following the considerations of Lemma~\ref{lem:8}, we obtain that $G^m (\om, x) \in R$, $m > n^{k+1}$, implies for $\om = (\om^1, \ldots, \om^{k+1})$
$$
( \om^{k+1}_{m-n^{k+1}} \ldots \om^{k+1}_{m-1} ) = ( 0 \ldots 0 ).
$$
This equation gives us the desired estimation $\e \le 2^{-n^{k+1}}$.

The part about the statistical attractor, as well as the part about perturbations, applies here without any changes. The proof of Theorem~\ref{thm:1} is completed. 

\section{Acknowledgements}

The first author is grateful to William Thurston, whose questions motivated writing of this paper.

The second author would like to thank Cornell University for hospitality and inspiring atmosphere which made possible the writing of this paper. 

%
\bigskip
Denis Volk (dire.ulf@gmail.com) \\
Institute for Information Transmission Problems \\
Russian Academy of Sciences \\
19 Bolshoy Karetny pereulok \\
127994 Moscow Russia
\\
\bigskip

\noindent Yulij Ilyashenko (yulij@math.cornell.edu) \\
Steklov Institute \\
8 Gubkina street \\
117966 Moscow Russia 

\end{document}